\documentclass{elsarticle}

\usepackage{amssymb,amsmath,amsfonts,amsthm}
\usepackage{amscd}
\usepackage{latexsym}
\usepackage{graphicx}
\usepackage{url}
\newenvironment{psmallmatrix}
  {\left(\begin{smallmatrix}}
  {\end{smallmatrix}\right)}

\newtheorem{theorem}{Theorem} 
\newtheorem{lemma}{Lemma} 
\newtheorem{proposition}{Proposition} 
\newtheorem{example}{Example} 
 
\newtheorem{corollary}{Corollary} 
 
\newtheorem{remark}{Remark}

\journal{Theoretical Computer Science}

\begin{document}

\begin{frontmatter}

\title{Morphic words and equidistributed sequences}


\author[I2M]{M\'elodie Andrieu}
\ead{melodie.andrieu@gmail.com}

\author[I2M]{Anna E. Frid}
\ead{anna.e.frid@gmail.com}
\ead[url]{http://iml.univ-mrs.fr/~frid/}

\address[I2M]{Aix Marseille Univ, CNRS, Centrale Marseille, I2M, Marseille, France}

\begin{abstract}
The problem we consider is the following: Given an infinite word $w$ on an ordered alphabet, construct the sequence $\nu_w=(\nu[n])_n$, equidistributed on $[0,1]$ and such that $\nu[m]<\nu[n]$ if and only if $\sigma^m(w)<\sigma^n(w)$, where $\sigma$ is the shift operation, erasing the first symbol of $w$. The sequence $\nu_w$ exists and is unique for every word with well-defined positive uniform frequencies of every factor, or, in dynamical terms, for every element of a uniquely ergodic subshift. In this paper we describe the construction of $\nu_w$ for the case when 
the subshift of 
$w$ is generated by a morphism of a special kind; then we overcome some technical difficulties to extend the result to all binary morphisms. The sequence $\nu_w$ in this case is also constructed with a morphism.

At last, we introduce a software tool which, given a binary morphism $\varphi$, computes the morphism on extended intervals and first elements of the equidistributed sequences associated with fixed points of $\varphi$.
\end{abstract}

\end{frontmatter}
\section{Introduction}
Consider an infinite word $w=w[0]w[1]\cdots w[n] \cdots$ on an ordered alphabet $\Sigma$;  here $w[n]\in \Sigma$. Suppose that the uniform frequency $\mu(u)$ of every factor $u$ of $w$ exists and is strictly positive, that is, that the dynamical system (subshift) associated with $w$ is uniquely ergodic, and $\mu$ is the unique ergodic measure on it (see \cite{fer_mont} for a discussion of this notion). Now for a factor $u$ of $w$, define $\nu(u)$ as the sum of measures $\mu$ of all words of the same length as $u$ which are lexicographically less than $u$, and $\nu(w)$ as the limit $\nu(w)=\lim_{n \to \infty} \nu (w[0]\cdots w[n])$.

The function $\nu$ on infinite words has been considered by Lopez and Narbel \cite{ln} from the dynamical point of view. 
On the other hand, as 
it was proven in \cite{words2015},
for every appropriate word $w$, the sequence $(\nu(\sigma^n(w)))_{n=0}^{\infty}$, where $\sigma$ is the shift operation, is uniformly distributed on $[0,1]$, and moreover, for $n\neq m$, we have $\nu(\sigma^n(w))\neq \nu(\sigma^m(w))$. This makes it possible to call the sequence $(\nu(\sigma^n(w)))_{n=0}^{\infty}$ the canonical representative of the {\it infinite permutation} defined by the shifts of $w$. Infinite permutations in this sense were introduced in \cite{ff}; as for permutations defined by words, their study was initiated independently by Makarov \cite{mak1, mak_tm, mak_st} and by Bandt, Keller and Pompe \cite{bkp}; see also \cite{elizalde} and the monograph \cite{a} summarizing that approach. The fact that the sequence $\{\nu(\sigma^n(w))\}_{n=0}^{\infty}$ is uniformly distributed means in particular that the respective permutation is also equidistributed (see \cite{dlt_subm} for the definition and discussion of an equidistributed permutation).

In this paper, given a morphism $\varphi$ with several nice properties, we describe how to find $\nu(w)$ for any infinite word $w$ from the respective subshift $L_{\varphi}$, and in particular for a fixed point $w_{\varphi}$ of $\varphi$. This result generalizes the Makarov's construction for the Thue-Morse word \cite{mak_tm}. A previous result in this direction, stated in combinatorial terms and considering not the whole subshift but just a fixed point of the morphism, can be found in \cite{words2015}. 

The next result of the paper concerns the binary case: if the morphism is binary, even if it does not belong to the ``nice'' class, our technique can be adapted to it. We also support the binary case by a software tool.

After introducing usual definitions in Section \ref{s:def0} and the object of our study in Section \ref{s:def}, we have to devote Section \ref{s:morph} to a discussion of properties of morphic subshifts. Section \ref{s:ext} starts with a correct extension of the interval $[0,1]$ to a wider set, which is needed to distinguish images of consecutive elements of the subshift. It also contains the first of main results of the paper, Theorem \ref{t:main}, 
giving a way 
to construct a morphism on numbers corresponding to the initial morphism $\varphi$, and a sequence $\nu_w$ for any element $w \in L_{\varphi}$. The construction is supported by examples.

Section \ref{s:kreg} contains a discussion of the case when the morphism $\varphi$ is $k$-uniform and thus its fixed point $w$ is $k$-automatic \cite{a_sh}. It is proved that in this case, the sequence $\nu_w$ is $k$-regular (see \cite{a_sh_reg} for the respective definitions).

The construction from Theorem \ref{t:main} works only for a restricted class of morphisms. However, in Section \ref{s:binary} we use some additional machinery to extend this result to any binary morphism. So, given a binary morphism, we know how to construct an equidistributed sequence corresponding to its fixed point(s) and to any element of the respective subshift.

At last, in Section \ref{s:experiment}, we discuss and refer to a software tool developed to compute the morphisms on numbers and sequences described in the paper. 

\section{Definitions and notation}\label{s:def0}

We consider infinite words $w=w[0]w[1]\cdots w[n] \cdots$, where $w[i]\in\Sigma$, on an ordered alphabet $\Sigma$. In this paper, we usually take $\Sigma=\{a,b,c,\ldots\}$, under the convention that $a<b<c<\cdots$. The order of symbols of $\Sigma$ naturally extends to the lexicographic order on finite and infinite words.

The factor $w[i]\cdots w[j]$ of a finite or infinite word $w$, where $j\geq i$, is denoted also by $w[i..j]$. The set of all factors of $w$ of length $n$ is denoted by Fac$_n(w)$.

The set of infinite words over $\Sigma$ is denoted by $\Sigma^{\omega}$.  As usual, the shift operation $\sigma$ corresponds to erasing the first symbol: $\sigma(w[0]w[1]\cdots w[n]\cdots)=w[1]w[2]\cdots w[n+1] \cdots$. Given an infinite word $w\in \Sigma^{\omega}$, we denote by $L_{w}$ the closure of the orbit of $w$ under $\sigma$. The dynamical system $(L_{w},\sigma)$ is called a {\it subshift} generated by $w$. 

An infinite word $w$ and its subshift $L_w$ are called {\it ultimately periodic} if $w=uvvvv\cdots$ for some finite words $u$ and $v$. If a word (or subshift) is not ultimately periodic, it is called {\it aperiodic}.

Given a set $S$, we denote by $S^*$ the set of finite concatenations of elements of $S$. In particular, if $S$ is an alphabet, $S^*$ is the set of finite words on $S$; but if $S$ is an interval of reals, $S^*$ is the set of finite sequences of numbers from $I$. A {\it morphism} $f:S^* \to S^*$ is a mapping which preserves concatenation, so that $f(xy)=f(x)f(y)$ for all $x,y \in S$. Clearly, a morphism is defined by its values on all elements on $S$.

Consider a morphism $\varphi: \Sigma^* \mapsto \Sigma^*$, where $\Sigma$ is an alphabet. If the image of a symbol of $x\in\Sigma$ starts with $x$, the morphism $\varphi$ admits a finite or right infinite fixed point $w_x$ starting with $x$ and defined as the limit $\lim_{n \to \infty} \varphi^n(x)$. If in addition $\varphi(x)\neq x$, $\varphi$ has no other fixed points starting with $x$. 

If the fixed point $w_x=\varphi(w_x)=\lim_{n \to \infty} \varphi^n(x)$ is infinite, it is called also a {\it pure morphic} infinite word, and the associated subshift $(L_{w_x},\sigma)$ is called a {\it pure morphic} subshift. In most cases considered in this paper (in particular, when the morphism is {\it primitive}, see the definition in Section \ref{s:morph}), the subshift does not depend on the letter $x$ and the fixed point $w_x$ but is uniquely defined by $\varphi$. In this case, it is denoted by $(L_{\varphi},\sigma)$, and its set of factors of length $n$ is denoted by Fac$_n(L_{\varphi})$.

Note that every element $u$ of a pure morphic subshift $(L_{\varphi},\sigma)$ can be obtained from the $\varphi$-image of another element $v=v[0]v[1]\cdots \in L_{\varphi}$ by the shift operation applied $p$ times, where $p \geq 0$. Moreover, we can choose $p$ to be less than the length of $\varphi(v[0])$: here of course we suppose that $\varphi(v[0])$ is not empty. So, $u=\sigma^p(\varphi(v))$, where $0\leq p < |\varphi(v[0])|$. Note that in the general case, $v$ and $p$, as well as $v[0]$ for a given $p$, are not unique.

A word $w$ is called {\it recurrent} if each of its factors $s=w[i..j]$ appears in it an infinite number of times. If in addition the distances between successive occurrences of $s$ are bounded, the word is called {\it uniformly recurrent}. As it is well-known, an infinite word $w$ is uniformly recurrent if and only if the associated subshift $(L_w,\sigma)$ is {\it minimal}, meaning that $L_w$ does not contain any proper subset which would be closed under $\sigma$. An even stronger condition is the existence of the unique $\sigma$-invariant probability measure $\mu$ on $L_w$, which is equivalent to the existence of uniform positive frequencies of all factors. In this case, the word $w$ and the dynamical system $(L_w,\sigma)$ generated by it are called {\it uniquely ergodic}. 

\section{Equidistributed sequences arising from words}\label{s:def}

Note that an infinite word $w$ is ultimately periodic if and only if $\sigma^m(w)=\sigma^n(w)$ for some $m\neq n$; so, if $w$ is aperiodic, for each $m\neq n$ we have either $\sigma^m(w)>\sigma^n(w)$ or $\sigma^m(w)<\sigma^n(w)$.

Consider an aperiodic word $w$ with well-defined non-zero uniform frequency $\mu(u)$ of every factor $u$. The subshift $L_w$ is uniquely ergodic, and its unique ergodic measure $\mu$ is completely determined by the frequencies $\mu(u)$ which can be interpreted as the values of the measure on cylinders: here a cylinder $[u]$ is the set of elements of $L_w$ starting with a word $u$. 

Now let us associate with an infinite word $v\in L_w$, that is, with an infinite word with the same set of factors that $w$, the measure
\[\nu(v)=\mu([w_{min},v]),\]
of the interval $[w_{min},v]$: here $w_{min}$ is the lexicographically minimal element of the subshift $L_w$, existing since the set $L_w$ is closed, and the interval $[w_{min},v]$ is  defined as the set of all infinite words from $L_w$ which are greater than or equal to $w_{min}$ and less than or equal to $v$.

The mapping $\nu: L_w \mapsto [0,1]$ is well-defined, and moreover, since among the shift images of $w$ the frequency of words from the interval $[w_{min},v]$ is the same as in the whole subshift,
\[\nu(v)=\lim_{n\to\infty} \frac{\#\{k|\sigma^k(w)<v,k\leq n\}}{n}.\]
Comparing it to the definition of a uniformly distributed (or, which is the same, equidistributed) sequence $(x[n])$ on an interval $[a,b]$, meaning that for every interval $[c,d]\subseteq [a,b]$, we have
\[\lim_{n \to \infty} \frac{\#(\{x[0],\cdots,x[n]\}\cap [c,d])}{n}=\frac{d-c}{b-a},\]
we see that for all $v \in L_w$, the sequence $\nu_v=(\nu(\sigma^n(v)))_{n=0}^{+\infty}$ is equidistributed on $[0,1]$. Indeed, since our subshift is uniquely ergodic, the proportion of words which are less than or equal to $\sigma^n(v)$ is the same in $v$, $w$ and the shift $L_w$ in total; see also a discussion in \cite{words2015}.

We will denote the real number $\nu(\sigma^n(w))$ by $\nu[n]$. By the construction, the sequence $\nu_w=(\nu[0],\nu[1],\ldots)$ is unique for every uniquely ergodic infinite word $w$.

The mapping $\nu$ has been considered by Lopez and Narbel in \cite{ln}. On the other hand, the equidistributed sequence $\nu_w$ for a sequence $w$, under another notation, was considered in \cite{dlt_subm} because of its relation to the {\it infinite permutation} generated by $w$; see \cite{mak1,mak_st,elizalde} for a discussion of infinite permutations defined by words. 

\begin{example}\label{e:tm}{\rm
The famous Thue-Morse word $w_{tm}=abbabaabbaababba\cdots$ \cite{a_sh_tm} is defined as the fixed point starting with $a$ of the morphism $\varphi_{tm}: a \to ab, b \to ba$. The sequence $\nu_{tm}=\nu_{w_{tm}}$ is equal to the fixed point
\[\frac{1}{2},1,\frac{3}{4},\frac{1}{4},\frac{5}{8},\frac{1}{8},\frac{3}{8},\frac{7}{8},\cdots,\]
of the  morphism $f_{tm}: [0,1]^*
\mapsto [0,1]^*$:
\begin{equation}\label{eq:tm}
 f_{tm}(x)=\begin{cases}
                                                          \frac{x}{2}+\frac{1}{4}, \frac{x}{2}+\frac{3}{4}, \mbox{~if~} 0 \leq x \leq \frac{1}{2},\\
                             \frac{x}{2}+\frac{1}{4}, \frac{x}{2}-\frac{1}{4}, \mbox{~if~} \frac{1}{2} < x \leq 1.
                                                         \end{cases}
\end{equation}
Note that morphisms on intervals have been discussed in Section \ref{s:def0} and, as any other morphisms considered in this paper, they transform a concatenation into a concatenation.

This construction (or, more precisely, a similar construction on the interval $[-1,1]$) was found by Makarov in 
2009 \cite{mak_tm}; below in Section \ref{s:ext} we shall prove its correctness as a corollary of a more general statement, Theorem \ref{t:main}.

In particular, we see from that construction that $\nu(0)=1/2$ and $\nu(1)=1$, which means that the Thue-Morse word is the maximal element of its subshift starting with $a$, and if we erase the first symbol from it, the result is the lexicographically maximal element of the subshift $L_{tm}$. These are known results (see, e.~g., \cite{berstel,johnson}).}
\end{example}

Note also that due to the symmetry between $a$ in $b$, the value of 
$\nu(w'_{tm})$ of the other fixed point $w'_{tm}=baababbaabba\cdots$ of the same morphism $\varphi_{tm}$ is also $1/2$. So, the mapping $\nu$ is not injective on the respective subshift $L_{tm}$. As it was discussed in \cite{ln}, this is a typical situation and 
it can be resolved by extending $[0,1]$ to a new wider domain described later in Section \ref{s:ext}. However, if we consider just a recurrent infinite word $w$ and its orbit, that is, the set of its shifts, and not the whole dynamical system to which it belongs, it is not necessary. Indeed, two infinite words with the same value of $\nu$ can never appear in the same orbit due to the following statement.

\begin{proposition}\label{c:maxmin} {\rm \cite{words2015}} 
Let $w$ be a recurrent aperiodic word and $u$ and $v$ be two of its
factors. Then the orbit of $w$ cannot contain at the same time the
lexicographically maximal word from $L_w$ starting with $u$ and
 the lexicographically minimal word from $L_w$ starting
with $v$.
\end{proposition}

In this paper, we consider only uniformly recurrent and, moreover, uniquely ergodic words, so, Proposition \ref{c:maxmin} can always be used.

\section{Properties of morphic symbolic subshifts}\label{s:morph}
In this section, we define the class of morphisms such that we can directly generalize the Thue-Morse construction above to their subshifts: namely, these are {\it primitive order-preserving} morphisms with {\it separable} subshifts. For such a morphism, we construct an interval
morphism similar to the Thue-Morse construction from Example \ref{e:tm}, and prove its correctness. The considered family of morphisms includes in particular all morphisms considered by Valyuzhenich
\cite{val}, and much more. Note that similar definitions have been introduced in \cite{words2015}, but here they are updated to better fit our wider goals.

Consider an alphabet $\Sigma=\{a_1,\ldots,a_q\}$ and let $\varphi: \Sigma^*\mapsto \Sigma^*$ be a
morphism with an aperiodic fixed point $u=\varphi(u)$ 
starting with a letter $a$. 

The matrix $M$ of a morphism $\varphi$ on a $q$-letter alphabet is a $q\times
q$-matrix whose element $m_{ij}$ is equal to the number of
occurrences of $a_i$ in $\varphi(a_j)$. The matrix $M$ and the morphism
$\varphi$ are called {\it primitive} if in some power $M^n$ of $M$
all the entries are positive, i.e.,  
for every $b \in \Sigma$ all the symbols of
$\Sigma$ appear in $\varphi^n(b)$ for some
$n$. The classical Perron-Frobenius theorem says that every primitive
matrix has a dominant positive {\it Perron-Frobenius eigenvalue}
$\theta$ such that $\theta>|\lambda|$ for any other eigenvalue
$\lambda$ of $M$. It is also well-known \cite{queffelec}
that a fixed
point of a primitive morphism is uniquely ergodic; moreover, for every sequence of factors $(v_n)$ of $L_{\varphi}$ of increasing length, the limit
\[\lim_{n\to \infty} \frac{|\varphi(v_n)|}{|v_n|}\]
exists and is equal to $\theta$.

Note in particular that every primitive morphism is non-erasing, which means that the images of all symbols are non-empty.
\begin{example}\label{e:fib}
 {\rm The Thue-Morse morphism is primitive with the matrix $\begin{psmallmatrix}
                                                               1&1\\ 1&1
                                                              \end{psmallmatrix}$. 
The Fibonacci morphism $\varphi_f: a \to ab, b\to a$ is primitive with the matrix 
$M=\begin{psmallmatrix}
1&1\\ 1&0
 \end{psmallmatrix}$: $M$ is not positive, but $M^2=\begin{psmallmatrix}
2&1\\ 1&1
 \end{psmallmatrix}$ is. The Sierpinski morphism $a \to aba, b \to bbb$ is not primitive since in all the powers of its matrix $\begin{psmallmatrix} 2&1\\ 0&3 \end{psmallmatrix}$, the left lower element is 0, and indeed, $a$ never appears in images of $b$.
}
\end{example}

We say that a morphism $\varphi$ is {\it order-preserving on an infinite
word} $u$ if for any $n,m>0$ we have $\sigma^n(u)<\sigma^m(u)$ if and only
if $\varphi(\sigma^n(u))<\varphi(\sigma^m(u))$; here $<$ denotes the
lexicographic order. A morphism is called {\it order-preserving} if it is
order-preserving on all infinite words, or, equivalently, if for any
infinite words $u$ and $v$ we have $u<v$ if and only if
$\varphi(u)<\varphi(v)$. If this property holds only for $u,v \in L_{\varphi}$, we say that $\varphi$ is {\it order-preserving on its subshift}. Note that in \cite{words2015}, order-preserving morphisms were called {\it monotone}.
\begin{example}\label{e:mon}
 {\rm The Thue-Morse morphism $\varphi_{tm}$ is order-preserving since $ab=\varphi_{tm}(a)<\varphi_{tm}(b)=ba$.  The Fibonacci morphism from Example \ref{e:fib} is not order-preserving since
 $ab=\varphi_f(a)>\varphi_f(ba)=aab$, whereas $a<ba$. At the same time, $\varphi_f^2: a \to aba, b \to ab$
 is order-preserving since for all $x,y\in\{a,b\}$ we have $\varphi^2_f(ax)=abaax'<ababy'=\varphi_f(by)$, where $x',y'\in\{a,b\}^*$. So, to use our construction on the Fibonacci word $u_f=\varphi_f(u_f)=abaab\cdots$  we should consider $u_f$ as the fixed point of $\varphi_f^2$ which is order-preserving.}
\end{example}

The last condition on the morphism $\varphi$, or, more precisely, on the subshift $(L_\varphi,\sigma)$, is to be {\it
separable}. To define this property, consider an element $u \in L_\varphi$ and all the ways to represent it as $\sigma^p(\varphi(u'))$ with $u' \in L_\varphi$ and $0\leq p<|\varphi(a)|$, where $a =u'[0]$ is the first symbol of $u'$. At least one such pair $(u',p)$ exists by the definition of $L_\varphi$. If this pair is unique, we call the pair $(a,p)$ the {\it type} $\tau(u)$ of $u$. A subshift is {\it typable} if for all elements $u \in L_{\varphi}$, the type of $u$ exists. If in addition the $\sum_{a\in \Sigma}|\varphi(a)|$ possible types can be ordered so that for all $u,v \in L_\varphi$ with $\tau(u)<\tau(v)$, we always have $u<v$, we say that the subshift $L_\varphi$ is {\it
separable}.

\begin{example}\label{tm_order}
 {\rm The Thue-Morse subshift $L_{tm}$ is separable. Indeed, first, any two consecutive $a$s (or $b$s) in its element determine a boundary between images of letters and thus all such boundaries. Also, the last symbols of $\varphi_{tm}(a)$ and $\varphi_{tm}(b)$ are different, the incomplete image of a symbol in the beginning can also be uniquely reconstructed, so, the morphism is typable. Moreover, if $\tau(u)=(a,0)$ and $\tau(v)=(b,1)$,
 we always have $u>v$, i.e., 
 all $a$s which are first symbols of $\varphi_{tm}(a)=ab$ give greater words than $a$s which
 are second symbols of $\varphi_{tm}(b)=ba$. The situation with $b$s is symmetric, so, we can order the types as
$(b,1)<(a,0)<(b,0)<(a,1)$ to have $u<v$ whenever $\tau(u)<\tau(v)$ for $u,v\in L_{tm}$.
}
\end{example}


\begin{example}\label{ex:nonsep}
 {\rm
The subshift $L$ generated by the morphism $\varphi: a \to aab, b \to abb$ is not typable because of the common suffix $b$ of images of letters. Indeed, consider 
a {\it special} infinite word $u$ such that $au,bu \in L$: such a word exists since the subshift is not periodic. Then the word $b \varphi(u)$ belongs both to $\sigma^2(\varphi(au))$ and to $\sigma^2(\varphi(bu))$, so that its type is not well-defined.
} 
\end{example}

\begin{example}\label{e:aabab}
 {\rm The subshift generated by the morphism $\varphi: a\to aabab, b\to bba$ is typable but not separable. Indeed, consider $u_1=abaa\cdots=\sigma^3(\varphi(aa\cdots))$, $u_2=ababaa\cdots=\sigma(\varphi(aa\cdots))$,  $u_3=abbb\cdots=\sigma^3(\varphi(ab\cdots))$. Then $u_1<u_2<u_3$ whereas $\tau(u_1)=\tau(u_3)=(a,3)$ and $\tau(u_2)=(a,1)$.}
\end{example}

For recurrence and, in some cases, precise formulas for the frequencies of factors in fixed points of morphisms, see \cite{queffelec,dol}.

In what follows, given a primitive order-preserving separable morphism $\varphi$ and the respective minimal subshift $L_{\varphi}$, we define a mapping which allows to build the sequence $\nu_w$, and in particular its first value $\nu(w)$, for any infinite word $w \in L_{\varphi}$. However, to do it, we first have to consider the extended domain to make the mapping $\nu$ injective.

\section{Extended intervals and morphisms}\label{s:ext}
As Lopez and Narbel showed in \cite{ln}, and as we discussed above just after Example \ref{e:tm}, the mapping $\nu: L\mapsto [0,1]$ defined in Section \ref{s:def} for any minimal subshift $L$ is surjective but not injective. In the Thue-Morse example, the image of the greatest word starting from $a$, which is $w_{tm}$ itself, is $1/2$, as well as the image of the smallest word starting from $b$. As it was proved in \cite{ln}, this happens exactly with {\it consecutive} words, or, equivalently, for consecutive cylinders. Recall that for a finite word $u$, a cylinder $[u]$ here is the set of all infinite words from $L$ starting from $u$. Finite words $u_1, u_2$ (and their cylinders $[u_1]$, $[u_2]$) are called {\it consecutive} if $u_1<u_2$ and there is no word $w \in L$ such that $w_1<w<w_2$, where $w_1$ is the greatest element of $L$ starting with $u_1$ and $w_2$ is the smallest element of $L$ starting with $u_2$. The infinite words $w_1$ and $w_2$ are also called consecutive. As it was proved in \cite{ln}, 
consecutive 
infinite words are exactly words $w_1\neq w_2$ for which $\nu(w_1)=\nu(w_2)$. Every pair of infinite consecutive words corresponds to a pair (or, more precisely, a countable number of pairs) of consecutive cylinders. For example, in the Thue-Morse subshift, the words $w_{tm}$ and $w'_{tm}$ are consecutive, as well as the respective cylinders $[a]$ and $[b]$, or $[abb]$ and $[baab]$, or any other pair of cylinders corresponding to prefixes of respectively $w_{tm}$ and $w'_{tm}$. 


Let $Z$ be the set of $\nu$-images of elements of consecutive pairs of words: it is a countable set since a consecutive pair can be defined by two finite (consecutive) words. 
To make the mapping $\nu$ injective and following \cite{ln}, we replace $Z$ in $[0,1]$ by two copies $Z^-$ and $Z^+$, and thus consider $\nu$ as a mapping from $L$ to the associated extended interval $X=X_L=([0,1]\backslash Z)\cup Z^- \cup Z^+$. Here for each pair $w_1<w_2$ of consecutive words with $\nu(w_1)=\nu(w_2)=x$ we denote $\nu(w_1)=a^-\in Z^-$ and $\nu(w_2)=a^+ \in Z^+$. It is natural to set $a^-<a^+$ and to make them both inherit from $[0,1]$ the relation with other elements of $X$. To unify the notation, we may also say for a number $a\in [0,1]\backslash Z$ that $a^-=a^+=a$.

Let $\varphi$ be a primitive order-preserving morphism  on an ordered alphabet $\Sigma=\{a_1,\ldots,a_q\}$, $a_1<\cdots<a_q$, with a separable subshift $(L,\sigma)$, and $X_L$ be the associated extended interval defined above. We will define a morphism on $X_L$ corresponding to $\varphi$, thus extending to $X_L$ a construction from \cite{words2015}.

Denote by $\mu=(\mu_1,\ldots,\mu_{q})$ the vector of measures of cylinders $[a_i]$ in $L$, or, which is the same, of
frequencies of symbols in any element of $L$. Since the morphism is primitive, these measures exist and are not equal to 0. 
Denote the intervals $I_{a_1}=[0,\mu_1^-]$,
$I_{a_2}=[\mu_1^+,(\mu_1+\mu_2)^-]$, $\ldots$, $I_{a_q}=[(1-\mu_{q})^+,1]$, $I_a \subset X_L$.

Now let us take all the $k=\sum_{i=1}^{q} |\varphi(a_i)|$ types of
elements of $L$ and denote them according to their order:
\[\tau_1< \tau_2 < \cdots < \tau_k,\]
with $\tau_i=(b_i,p_i)$. Types and their order exist since the subshift is separable.

For each $\tau_i=(b_i,p_i)$, the frequency of factors of type $\tau_i$ in the subshift is equal to $l_i=\mu_{b_i}/\theta$, where $\theta$ is the Perron-Frobenius eigenvalue of $\varphi$. Indeed, given a word $u$ from the subshift $L_{\varphi}$, consider its $\varphi$-image $v$ interpreted as a word on the alphabet $\varphi(\Sigma)$. By the construction, occurrences of $\varphi(b_i)$ to $v$ correspond to occurrences of $b_i$ to $u$. But if now we interpret $v$ as a word on $\Sigma$, $v \in L_{\varphi}$, we see that its prefix corresponding to the prefix of $u$ of length $n$ has a length which grows as $\theta n$ with $n\to \infty$. So, factors of type $\tau_i$ occur in $v$ exactly $\theta$ times rarer than $b_i$ in $u$. Since the frequencies do not depend on the choice of an element of $L_{\varphi}$, we get the formula $l_i=\mu_{b_i}/\theta$.

Denote
\[J_1=[0,l_1^-], J_2=[l_1^+,(l_1+l_2)^-],\ldots, J_k=[(1-l_k)^+,1];\]
so that in general, $J_i=[(\sum_{m=1}^{i-1}l_m)^+,(\sum_{m=1}^{i}l_m)^-]$.
We will also denote $J_i=J_{b_i,p_i}$. 

The interval $J_i$ is the range of $\nu(u)$ corresponding to elements of $u \in L$  of type $\tau_i$. Note that the first symbol of such a word $u$ is always the symbol number $p_i+1$ of $\varphi(b_i)$ (the range of $p_i$ for a given $b_i$ is from 0 to $|\varphi(b_i)|-1$). So, the union of elements $J_i$ corresponding to this element equal to $a_m$ is exactly $I_m$ for every $m$. In particular, each $J_i$
is a subinterval of some $I_{m}$. By the construction, all the ends of these intervals are in $Z^-$ or $Z^+$, and thus the intervals $J_i$ do not intersect: the ends $a^-$ and $a^+$ of consecutive intervals correspond to consecutive words from $L$.

Now we define the morphism $f: X_L^* \mapsto X_L^*$ as
follows: For $x \in I_a$ we have
\[f(x) = f_{a,0}(x),\ldots f_{a,|\varphi(a)|-1}(x).\]
Here $f_{a,p}$ is the increasing 
affine bijection $f_{a,p}:I_a \mapsto
J_{a,p}$: If $I_a=[x_1^+,x_2^-]$ and $J_{a,p}=[y_1^+,y_2^-]$, then
\begin{equation}\label{e:psi}
f_{a,p}(x)=\frac{y_2-y_1}{x_2-x_1}(x-x_1)+y_1.
\end{equation}
Here, by the convention, the image of any $x \in Z^-$ ($x \in Z^+$) is $(f_{a,p}(x))^-$ (respectively, $(f_{a,p}(x))^+)$. Note that the slope $\displaystyle \frac{y_2-y_1}{x_2-x_1}$ of the affine mapping $f_{a,p}(x)$ is equal to $1/\theta$ since the interval $J_{a,p}$ is $\theta$ times shorter than $I_a$.

The meaning of intervals $I_a$ and of the morphism $f$ is explained in the following proposition following directly from the construction.
\begin{proposition}
 Let $d: X_L^* \mapsto \Sigma^*$ be the morphism defined by $d(x)=a$ whenever $x \in I_a$. Then for all $x \in X_L$ we have 
$d(f(x))=\varphi(d(x))$.
\end{proposition}

This proposition means in particular that the lengths of $\varphi$-images of letters and of $f$-images of reals from the respective intervals are synchronized. So, the following statement holds.

\begin{proposition}\label{p:synch}
Given a word $w \in L_{\varphi}$, where $\varphi$ is primitive and order-preserving and $L_\varphi$ is separable, the following statements are equivalent: \begin{itemize}
    \item 
    The letter $\varphi(w)[n]$ is the letter indexed $p$ of the $\varphi$-image of the letter indexed $n'$ of $w$, and
    \item 
    The number $f(\nu_w)[n]$ is the number indexed $p$ of the $f$-image of the number indexed $n'$ of $\nu_w$.
\end{itemize}

\end{proposition}

\begin{example}
 {\rm 
The Thue-Morse morphism on $[0,1]$ from Example \ref{e:tm} can now be more correctly redefined on the respective set $X_{tm}$. Here $1/2$, which is the frequency of $a$, is one of the numbers which is doubled, as well as all binary rationals from $(0,1)$. So, we have $I_a=[0,1/2^-]$, $I_b=[1/2^+,1]$, $J_{a,0}=[1/4^+,1/2^-]$, $J_{a,1}=[3/4^+,1]$,
$J_{b,0}=[1/2^+,3/4^-]$, $J_{b,1}=[0,1/4^-]$, and \eqref{eq:tm} can now be rewritten as
\begin{equation}\label{e:tm2}
 f_{tm}(x)=\begin{cases}
   f_{a,0}(x),f_{a,1}(x) \mbox{~for~} x \in I_a,\\
   f_{b,0}(x),f_{b,1}(x)  \mbox{~for~} x \in I_b,
                                                         \end{cases}
\end{equation}
where the used linear mappings on $X_{tm}$ are defined by \eqref{e:psi} and, of course, coincide on $[0,1]$ with those from \eqref{eq:tm}.
}
\end{example}
For an example concerning the Fibonacci word as the fixed point of the square morphism $\varphi^2_f$ from Example \ref{e:mon}, see \cite{words2015}; the only difference in the presentation should be extended intervals.

The following statements is one of the main results of the paper.

\begin{theorem}\label{t:main}
Let $\varphi$ be a primitive morphism defining a separable subshift $(L,\sigma)$ and order-preserving on it, $f$ be the morphism on extended intervals associated with $\varphi$ (and $L$), and $\nu_w \in X_L^{\omega}$ be the equidistributed sequence corresponding to a sequence $w \in L$. Then 
$f(\nu_w)=\nu_{\varphi(w)}$ (see the commutative diagram below).
\end{theorem}
\begin{center}
$\begin{CD}
w @>\varphi>> \varphi(w)\\
@VV\nu V @VV\nu V\\
\nu_w @>f>> \nu_{\varphi(w)}
\end{CD}$
\end{center}

{\sc Proof.}
We shall prove first, that the two sequences have the same order among elements, and second, that $f(\nu_w)$ is equidistributed on $[0,1]$. Since the equidistributed sequence on $[0,1]$ corresponding to any ordering of elements is at most unique (if it exists, each of its element is uniquely defined as the fraction of elements in the ordering smaller than it, see also a discussion after Definition 2.3 in \cite{dlt_subm}), this is sufficient. 

Suppose that $f(\nu_w)[n]<f(\nu_w)[m]$ for some $n,m \in \mathbb N$; our goal is to prove that $\nu_{\varphi(w)}[n]<\nu_{\varphi(w)}[m]$. 

Suppose first that $f(\nu_w)[n]$ and $f(\nu_w)[m]$ are situated in the same interval $J_{c,p}$. Since all such intervals are disjoint, this means that $f(\nu_w)[n]$ is obtained as $f_{c,p}(\nu_w[n'])$ and $f(\nu_w)[m]$ is obtained as $f_{c,p}(\nu_w[m'])$ for some $n',m' \in \mathbb N$, where $\nu_w[n']$, $\nu_w[m'] \in I_c$. So, by the definition of $I_c$, we have $w[n']=w[m']=c$. Moreover, since $f$ and $\varphi$ are synchronized as described in Proposition \ref{p:synch}, the symbols of $\varphi$-images of $w[n']$ and $w[m']$ numbered $p$ are $\varphi(w)[n]$ and $\varphi(w)[m]$. So, the inequality $f(\nu_w)[n]<f(\nu_w)[m]$ implies that $\nu_w[n']<\nu_w[m']$ since $f_{c,p}$ is an affine mapping with a positive slope, then that $\sigma^{n'}(w)<\sigma^{m'}(w)$ by the definition of $\nu_w$, then that $\varphi(\sigma^{n'}(w))<\varphi(\sigma^{m'}(w))$ since $\varphi$ is order-preserving on $L$, and $\sigma^p(\varphi(\sigma^{n'}(w)))<\sigma^p(\varphi(\sigma^{m'}(w)))$ since first $p$ symbols of $\varphi(\sigma^{n'}(w))$ and $\varphi(\sigma^{m'}(w))$ are equal (to the first $p$ symbols of $\varphi(c)$). But since $\varphi$ and $f$ are synchronized by Proposition \ref{p:synch}, $\sigma^p(\varphi(\sigma^{n'}(w)))=\sigma^n(\varphi(w))$ and $\sigma^p(\varphi(\sigma^{m'}(w)))=\sigma^m(\varphi(w))$. So, $\sigma^n(\varphi(w))< \sigma^m(\varphi(w))$, and, by the definition of $\nu_{\varphi(w)}$, $\nu_{\varphi(w)}[n]<\nu_{\varphi(w)}[m]$. 

Now suppose that $f(\nu_w)[n]\in J_{c_n,p_n}$ and $f(\nu_w)[m]\in J_{c_m,p_m}$, where $J_{c_n,p_n} \neq J_{c_m,p_m}$. These intervals are disjoint and correspond to types $(c_n,p_n)<(c_m,p_m)$. So, $f(\nu_w)[n]$ is the number indexed $p_n$ in the $f$-image of some $\nu_w[n']\in I_{c_n}$, and $f(\nu_w)[m]$ is the number indexed $p_m$ in the $f$-image of some $\nu_w[m']\in I_{c_m}$. By the construction of intervals $I_c$, this means that $w[n']=c_n$ and $w[m']=c_m$, and moreover, due to Proposition \ref{p:synch}, $\varphi(w)[n]$ is the symbol indexed $p_n$ of the $\varphi$-image of $w[n']$ and $\varphi(w)[m]$ is the symbol indexed $p_m$ of the $\varphi$-image of $w[m']$. Since the morphism $\varphi$ is separable, the type of $\sigma^n(\varphi(w))$ is thus equal to $(c_n,p_n)$, and it is less than the type of $\sigma^m(\varphi(w))$ equal to $(c_m,p_m)$. So, $\sigma^n(\varphi(w))<\sigma^m(\varphi(w))$ and thus $\nu_{\varphi(w)}[n]<\nu_{\varphi(w)}[m]$, which was to be proved.

We have proved that the sequences $f(\nu_w)$ and $\nu_{\varphi(w)}$ have the same order among elements; it remains to prove that $f(\nu_w)$ is equidistributed on $X_L$. Indeed, let us consider any subinterval $I$ of length $l$ of some interval $I_c$. The frequency of elements of $\nu_w$ which are in $I$ is $l$ since $\nu_w$ is equidistributed. Due to the definition of $f$, the images of the interval $I$ are $f_{c,1}(I),\ldots,f_{c,|\varphi(c)|}(I)$. These are intervals from $X_L$ of length $l/\theta$ each, where $\theta$ is the Perron-Frobenius eigenvalue of $\varphi$. The frequency of elements of $f(\nu_w)$ from each of these intervals is also $l/\theta$, since $f$ and $\varphi$ are synchronized in length and since $\theta$ is the limit of the ratio $|\varphi(w[0]\cdots w[n])|/n$ with $n\to \infty$; we use also the fact that the intervals $J_{x,p}$ form a disjoint partition of $X_L$. On the other hand, every subinterval $J$ of some $J_{c,p}$, where $p=1,\ldots,|\varphi(c)|$ is the $f_{c,p}$-image of a respective subinterval $I$ of $I_c$ which is $\theta$ times longer than it. So, the frequency of elements from $J$ in $f(\nu_w)$ is equal to the length of $J$. This is true for all subintervals of $J_{c,p}$ and thus by union for all subintervals of $X_L$, meaning exactly that the sequence $f(\nu_w)$ is equidistributed on $X_L$. \hfill $\Box$


\begin{corollary}
 Let $w =w_a \in L_{\varphi}$ be the fixed point of $\varphi$ starting with $a$. Then $\nu_w$ is the unique fixed point of $f$ starting from a number from $I_a$, which is the fixed point of $f_{a,0}$.
\end{corollary}
{\sc Proof.}
 First of all, the fixed point $w_a$ of $\varphi$ starting with $a$ is unique since the morphism is primitive and thus $\varphi(a)\neq a$. We know from Theorem \ref{t:main} that $\nu_{\varphi(w_a)}=f(\nu_{w_a})$, but since $w_a=\varphi(w_a)$, here it means that $\nu_{w_a}=f(\nu_{w_a})$, and so $w_a$ is a fixed point of $f$ starting with a number from $I_a$ which is a fixed point of $f_{a,0}$, the first applied interval morphism. Since $f_{a,0}$ is an affine function with the slope $1/\theta<1$, this fixed point is unique. \hfill $\Box$

\begin{example}{\rm
The sequence $\nu_{tm}$, $\nu_{tm}[k]\in X_{tm}$, corresponding to the Thue-Morse word $w_{tm}$ starting with $a$, is the fixed point starting with $1/2^-$ of the morphism \eqref{e:tm2}: 
\[1/2^-,1,3/4^-,1/4^-,5/8^-,1/8^-,3/8^-,7/8^-,\cdots.\]
The other fixed point $w_{tm}'$ of $\varphi_{tm}$, starting with $b$, corresponds to the fixed point of the same morphism \eqref{e:tm2} starting with $1/2^+$:
\[1/2^+,0,1/4^+,3/4^+,3/8^+,7/8^+,5/8^+,1/8^+,\cdots.\]
Compared to Example \ref{e:tm}, we see that the fixed points here differ from the original definitions of $\nu_{tm}$ and $\nu_{tm'}$ by signs $-$ or $+$ added to numbers.
}
\end{example}

\begin{remark}
 {\rm
The set $Z$ existing in $X_{tm}$ in two copies, $Z^+$ and $Z^-$, contains not only binary rationals from the fixed points above. As another example, consider the number $1/6$. We have $1/6=f_{b,1}(f_{a,1}(1/6))$, and thus $1/6$ corresponds to words $au$ starting from $a$ and satisfying $au=\sigma(\varphi_{tm}(\sigma(\varphi_{tm}(au))))=\sigma^3(\varphi_{tm}^2(au))=a\varphi_{tm}^2(u)$. This equation has two solutions $aw_{tm}<aw'_{tm}$. So, we have $\nu(aw_{tm})=1/6^-$ and $\nu(aw'_{tm})=1/6^+$. Since $1/6$ is the frequency of $aa$ in the Thue-Morse word, we see that $aw_{tm}$ is the maximal element of $L_{tm}$ starting from $aa$ and $aw'_{tm}$ is the minimal element of $L_{tm}$ starting from $ab$.
}
\end{remark}

Other examples, starting with the Fibonacci morphism, can be treated with the software tool described below in Section \ref{s:experiment} and available online.

\section{Morphism $f$ and $k$-regular sequences}\label{s:kreg}
Let $\varphi$ be a primitive order-preserving morphism, $w$ be some of its fixed points, and $L$ be the separable subshift generated by $w$. Since the morphism is primitive, the subshift $L$ does not depend on the choice of the fixed point of $\varphi$ and is minimal.

Consider the morphism $f: X_L^* \mapsto X_L^*$ described above and the sequence $\nu_w$ which is its fixed point corresponding to $w$. As we have discussed above, each mapping $f_{a,p}(x)$ from the definition of $f$ is an affine mapping sending the interval $I_a$ of length $\mu_a$ to the interval $J_{a,p}$ of length $\mu_a/\theta$, where $\theta$ is the Perron-Frobenius eigenvalue of $\varphi$. So, in the definition \eqref{e:psi} of $f_{a,p}$, we have $(y_2-y_1)/(x_2-x_1)=1/\theta$. We get the following
\begin{corollary}
Under the conditions of Theorem \ref{t:main}, the mappings $f_{a,p}$ from the definition of the morphism $f$ are of the form
\begin{equation}\label{e:th}
f_{a,p}(x)=x/\theta +C_{a,p},
\end{equation}
where $\theta$ is the Perron-Frobenius eigenvalue of $\varphi$ and $C_{a,p}$ is a constant defined, in the notation of \eqref{e:psi}, by
\[C_{a,p}=y_1-x_1/\theta.\] 
\end{corollary}

This statement is particularly interesting when the morphism $\varphi$ is $k$-uniform, that is, the length of all $\varphi$-images of letters is the same and equal to $k\geq 2$. Since $\theta$ is the limit of the ratio $|\varphi(w[0..n])|/n$, here we have $\theta=k$. The word $w$ is $k$-automatic (for the definition and discussion on $k$-automatic words and $k$-regular sequences, the reader is referred to \cite{a_sh,a_sh_reg}). 

\begin{lemma}
 Under the conditions of Theorem \ref{t:main}, if the morphism $\varphi$ is $k$-uniform, then the sequence $\nu_w$ is $k$-regular.
\end{lemma}
{\sc Proof.} First of all, the morphism $f$ is also $k$-uniform. With the definition \eqref{e:th} of each mapping $f_{a,p}$, we can write the following expression for an element $\nu_w[kn+p]$ of the sequence $(\nu_w[n])_{n=0}^{\infty}$, where $n\geq 0$ and $p\in\{0,\ldots,k-1\}$:
\begin{eqnarray*}
\nu_w[kn+p]&=&\frac{1}{k}\nu_w[n]+C_{w[n],p}\\
&=&\frac{1}{k}\nu_w[n]+\sum_{q=0}^{k-1}\sum_{a \in \Sigma} C_{a,q}X(w[n]=a)X(q=p).
\end{eqnarray*}
Here $X(P)$ is the characteristic sequence of a property $P$, equal to 1 if the property holds and to 0 otherwise.

The sequence $w$ is $k$-automatic and thus $k$-regular, as well as the sequences $X(w[n]=a)$ and $X(w[n]=a)X(q=p)$. So, the sequence $\nu_w$ is also $k$-regular by Theorem 16.1.3 from \cite{a_sh_reg}. \hfill $\Box$

\section{The binary case}\label{s:binary}
In this section, we adapt Theorem \ref{t:main} to all binary morphisms with aperiodic uniformly recurrent fixed points. To do it, we first discuss what may happen. In this section, we consider the alphabet $\Sigma_2=\{a,b\}$ with $a<b$ and a morphism $\varphi:\Sigma_2^* \mapsto \Sigma_2^*$.

\subsection{Non-primitive case}
First of all, let us discuss the condition of the morphism $\varphi$ to be primitive. In fact, everywhere in the proof we used not the primitivity itself but the facts that the subshift $L_{\varphi}$ is uniquely ergodic, and for any finite factor $u$ of $L_{\varphi}$, the relation $|\varphi(u)|/|u|$ has a limit, denoted $\theta$ with $|u| \to \infty$. We also need the subshift to be aperiodic.

If $\varphi$ is primitive, all the conditions except for perhaps aperiodicity hold. Consider the case of non-primitive $\varphi$. It may have several fixed points with different orbit closures, but without loss of generality, suppose that $\varphi$ has an infinite fixed point $w$ starting with $\varphi(a)$ in which both letters appear. The condition that $\varphi$ is not primitive means then that $\varphi(b)\in b^*$.

\begin{proposition}\label{p:nonprim}
The fixed point $w$ and its orbit closure can be minimal and aperiodic only if $\varphi(b)=b$ and $\varphi(a)=axa$, where $x$ is a finite word containing $b$. In this case, $w$ is also uniquely ergodic, and there exists a limit
\begin{equation}\label{e:lim}
 \lim_{u\in {\rm Fac(w)}, |u|\to \infty }\frac{|\varphi(u)|}{|u|}=\theta.
\end{equation}
\end{proposition}

\noindent {\sc Proof.}
If $\varphi(b)$ is the empty word, then $w=\varphi(a)^{\omega}$ is periodic. If $\varphi(b)=b^k$ with $k \geq 2$, then $L_w$ is not minimal since contains $b^{\omega}$. The same is true if $\varphi(b)=b$ and $\varphi(a)$ ends with $b$. So, the only case when $L_w$ is minimal and can be aperiodic is $\varphi(b)=b$ and $\varphi(a)=axa$, where $x$ is a word containing $b$. In this case, $w$ is uniformly recurrent: indeed, the distance between two consecutive occurences of $a$ is bounded by $|x|+1$ and thus the distance between two consecutive occurrences of any factor of $\varphi^n(a)$ is bounded by $|\varphi^n(ax)|$. Moreover, due to Theorem 3 of \cite{durand}, $w$ is a primitive morphic sequence. In particular, it is uniquely ergodic, and the limit  \eqref{e:lim} exists. \hfill $\Box$

Note that in this case, the word $w$ and its subshift can be periodic if $\varphi(a)=(ab)^ka$ for some $k$. In other cases, we can work with $\varphi$, $w$ and $L_w$ exactly as if $\varphi$ were primitive.

\begin{example}\label{e:nonprim}{\rm
 If $\varphi(a)=aaba, \varphi(b)=b$, the minimal subshift $L_w$ generated by $w$ is the orbit closure of the fixed point $w$ of $\varphi$ starting with $a$: $$w=aabaaababaabaaabaaababaaba\cdots.$$
We can work with $L_w$ exactly as if $\varphi$ were primitive.
}
\end{example}
\subsection{Preserving order}
How can we work with a morphism which is not order-preserving? The following proposition shows that in the binary case, we can just replace it by its square, as we did in Example \ref{e:mon} for the Fibonacci morphism.
\begin{proposition}\label{p:monotoneSq}
For every binary morphism $\varphi$ with an aperiodic subshift, if $\varphi$ is not order-preserving, then $\varphi^2$ is.
\end{proposition}
\noindent {\sc Proof.} Borchert and Rampersad proved (see Theorem 15 in \cite{br}) that every aperiodic binary morphism is either order-preserving or order-reversing, the latter property meaning that $\varphi(u)>\varphi(v)$ whenever $u<v$ for infinite words $u,v$. Clearly, the square of an order-reversing morphism is order-preserving. \hfill $\Box$

So, if by chance a binary morphism is not order-preserving, its square is, and we can consider the subshift as generated by the square morphism. Now let us consider two different types of inseparability.

\subsection{Common suffixes}\label{ss:cs}
Here we consider the situation when the morphism $\varphi$ is not typable because of a non-empty common suffix of images of letters, like in Example \ref{ex:nonsep}. Since we consider only the minimal (or, which is the same, uniformly recurrent) aperiodic case, we can use Proposition \ref{p:nonprim} and see that if $\varphi$ is not primitive, the common suffix is empty. So, in our case, the morphism $\varphi$ is primitive. 

The following classical statement will be useful. We give its proof for the sake of completeness.

\begin{proposition}\label{p:pssp}
Suppose that $\varphi: a \to p_a s, b \to p_b s$, where $s$ is any common suffix of $\varphi(a)$ and $\varphi(b)$. Then $L_{\varphi}=L_{\varphi'}$, where $\varphi': a \to s p_a, b \to s p_b$.
\end{proposition}
\noindent {\sc Proof.} Clearly, both subshifts generated by $\varphi$ and $\varphi'$  contain the word $s$ and are closed under the operation sending a word $u=u_1 u_2\cdots u_n$ to the word $p_{u_1}s p_{u_2} s \cdots s p_{u_n}$. So, the intersection of the two sets of factors is infinite, and since both subshifts are minimal, they are equal. \hfill $\Box$

\begin{example}\label{ex:nonsep2}
{\rm For the morphism $\varphi: a \to aab, b \to abb$ from Example \ref{ex:nonsep}, we have $\varphi': a \to baa, b \to bab$.}
\end{example}

Clearly, if one of the images of symbols is not a suffix of the other, it is sufficient to apply Proposition \ref{p:pssp} once to get a new morphism, which is also primitive and with images of letters ending with different symbols, like in the previous example. If it is not the case, however, it can be necessary to apply the operation from Proposition \ref{p:pssp} several times:

\begin{example}\label{e:babab}
 {\rm 
If $\varphi: a \to ab, b\to babab$, then, to get a morphism with images of letters ending by different symbols, we follow three steps:
\begin{eqnarray*} \varphi: \begin{cases}
            a \to ab \\ b\to babab
           \end{cases}
&\to&
\varphi': \begin{cases}
            a \to ab \\ b\to abbab
           \end{cases} \\
&\to&
\varphi'': \begin{cases}
            a \to ab \\ b\to ababb
           \end{cases}
\to
\varphi''': \begin{cases}
            a \to ba \\ b\to babab.
           \end{cases} \end{eqnarray*}
}
\end{example}
Note that in terms used in \cite{richomme_conj}, $\varphi'$ is a {\it left conjugate} of $\varphi$; for conjugacy of morphisms see also \cite{seebold_conj}. The next proposition follows from results of \cite{richomme_conj}, but it does not take more space to give a new proof than to explain the relationship between used terminology. 
The meaning of the proposition is that the number of necessary steps is always finite.

\begin{proposition}\label{p:suffixes}
Let $\varphi$ be a primitive binary morphism with aperiodic subshift and, without loss of generality, $|\varphi(a)|\geq |\varphi(b)|$. Then after applying the operation from Proposition \ref{p:pssp} with the maximal common suffix at most $k$ times, where $k\leq \lfloor|\varphi(a)|/ |\varphi(b)| \rfloor +1$, we get a morphism whose images of letters end by different symbols.
\end{proposition}
\noindent {\sc Proof.} If one image of letter is not a suffix of the other, the statement is obvious and $k=1$ is enough. Suppose now the opposite: let $\varphi(a)=ps, \varphi(b)=s$. Suppose that we can continue to exchange the prefix and the suffix $s$ of the image of $a$ at least $(|\varphi(a)|+|\varphi(b)|)/ |\varphi(b)|$ times. It would mean that the word 
$sps=\varphi(ba)$ is periodic with the period $|s|$. But it is also periodic with the period $|ps|$, so, due to the Fine and Wilf theorem, it is periodic with the period $\gcd (|s|,|ps|)$. In particular, it means that $\varphi(a)=ps$ and $\varphi(b)=s$ are powers of the same word, and thus the subshift generated by $\varphi$ is periodic, a contradiction. So, if $k$ is the maximal number of replaced suffixes, we have $k<(|\varphi(a)|+|\varphi(b)|)/ |\varphi(b)|$, which means $k\leq \lfloor|\varphi(a)|/ |\varphi(b)| \rfloor +1$.
 $\hfill$ $\Box$

\begin{example}
 {\rm 
Note that nevertheless, successive transfers of the longest common suffixes to the left can touch more symbols than there are in the longer image of a letter. For example, consider $\varphi: a \to abbab, b\to bab$; then 
\begin{equation*} \varphi: \begin{cases}
            a \to abbab \\ b\to bab
           \end{cases}
\to
\varphi': \begin{cases}
            a \to babab \\ b\to bab
           \end{cases}
\to
\varphi'': \begin{cases}
            a \to babba \\ b\to bab.
           \end{cases}
\end{equation*}
Here the longer image of a letter is of length 5, and there are 3+3=6 letters replaced.
}
\end{example}
So, we may assume that given a primitive binary morphism $\varphi$, we can always transfer common suffixes of images of letters to the left until we get a morphism $\psi=\varphi^{(k)}$ with the same subshift and with images of letters ending by different symbols. To justify this passage completely, we should also describe how we can apply Theorem \ref{t:main} to $\varphi$ if we know how to do it for $\psi$. The following proposition gives a recipe for that.

\begin{proposition}
 Suppose that a binary morphism $\varphi$ is transformed to another morphism $\psi$ by a series of transfers of common suffixes to the left: $\varphi \to \varphi' \to \cdots \to \varphi^{(k)}=\psi$, where the suffix transferred at the passage from $\varphi^{(i)}$ to $\varphi^{(i+1)}$ is of length $p_{i+1}$. Then for every infinite word $w$, we have
$\psi(w)=\pi \varphi(w)$, where the word $\pi$ of length
$p=p_1+\cdots + p_k$ is the concatenation of all replaced common suffixes in order from right to left.
\end{proposition}
\noindent {\sc Proof.} For each step $i$, it is not difficult to see that 
$\varphi^{(i+1)}(w)=\pi_{i+1}\varphi^{(i)}(w)$, where $\pi_{i+1}$ is the replaced common suffix of length $p_{i+1}$. It remains to combine these arguments for all $i$ and to set $\pi=\pi_k \cdots \pi_2 \pi_1$. \hfill $\Box$

\begin{corollary}\label{c:opres}
The morphism $\varphi$ is order-preserving if and only if $\psi$ is order-preserving.
\end{corollary}

This corollary means just that we can successively take a square of our morphism if it is necessary to make it order-preserving, and then transfer common suffixes as we need.

\begin{example}\label{e:aababb}
 {\rm 
Consider the fixed point $w=aabaababb\cdots$ of the morphism $\varphi: a \to aab, b \to abb$ from Examples \ref{ex:nonsep} and \ref{ex:nonsep2}. To find the value $\nu(w)$ and all the sequence $\nu_w$, we pass to the morphism $\psi=\varphi': a \to baa, b \to bab$. The morphism $\psi$ falls into conditions of Theorem \ref{t:main}, and gives rise to the the following morphism on intervals:
\begin{eqnarray*}f(x)&=&\begin{cases}
           f_{a,0}(x),f_{a,1}(x),f_{a,2}(x) \mbox{~for~} x \in [0,1/2^-]\\
           f_{b,0}(x), f_{b,1}(x),f_{b,2}(x) \mbox{~for~} x \in [1/2^+,1]
          \end{cases}\\
&=&\begin{cases}
           x/3+1/2,x/3,x/3+1/6 \mbox{~for~} x \in [0,1/2^-]\\
           x/3+1/2, x/3+1/6,x/3+2/3 \mbox{~for~} x \in [1/2^+,1].
          \end{cases}
\end{eqnarray*}
The only fixed point $v=babbaabab\cdots$ of $\psi$ corresponds to the only fixed point $\nu_v$ of $f$ starting with the fixed point $\nu(v)=3/4$ of $f_{b,0}(x)=x/3+1/2$: $\nu_v=3/4,5/12,11/12,23/36,\ldots$. At the same time, the fixed point $w$ of $\varphi$ due to the previous proposition satisfies $w=\varphi(w)=\sigma(\psi(w))$, and due to Proposition \ref{p:synch} and Theorem \ref{t:main}, corresponds to the fixed point $\nu_w$ of $\sigma(f)$: $\nu_w=\sigma(f(\nu_w))$. In particular, it starts with the fixed point $\nu(w)=0$ of $f_{a,1}(x)=x/3$: $\nu_w=0,1/6^+,10/18^+,1/18^+,2/9^+,\cdots$. Here we have to add pluses to values since these are lower ends of intervals.
}
\end{example}

\begin{example}\label{e:lettertransf}{\rm
 Consider the morphism $\varphi: a \to ab, b \to babab$ from Example \ref{e:babab}. After moving $p=2+2+1=5$ symbols from right to left, we get an order-preserving separable morphism $\psi=\varphi''': a \to ba, b \to babab$ inducing the same subshift $L$.

Note that by the construction, for every $u \in L$, we have $\varphi(u)=\sigma^5(\psi(u))$. In particular, it is true for both fixed points $w_a$ and $w_b$ of $\varphi$. Let us start with the fixed point $w=w_a$ starting with $a$:
\[w=ab.babab.babab.ab.babab.ab.babab.babab.\ldots=\varphi(w)=\sigma^5(\psi(w))\]
(dots are put between $\varphi$-images of symbols for readability).
Since $w$ starts with $a$ and $\psi(a)$ is of length 2, we have
\[w=\sigma^5(\psi(w))=\sigma^3(\psi(\sigma(w))).\]
Passing to the sequences $\nu$, we see that
\[\nu_w=\sigma^3(f(\sigma(\nu_w))),\]
where $f$ is the morphism on extended intervals corresponding to $\psi$.
Denote $\nu_w=\nu[1]\nu[2]\nu[3]\cdots$. Here $\nu[1]\in I_a$, $\nu[2],\nu[3]\in I_b$ and so on along the word $w$. Then $\sigma(\nu_w))=\nu[2]\nu[3]\nu[4]\cdots$,  and 
\[f(\sigma(\nu_w))=f_{b,0}(\nu[2])f_{b,1}(\nu[2])\cdots f_{b,4}(\nu[2]) f_{b,0}(\nu[3])f_{b,1}(\nu[3])\cdots\]
At last, applying $\sigma^3$, we see that
\[\nu_w=\nu[1]\nu[2]\nu[3]\nu[4]\cdots=f_{b,3}(\nu[2])f_{b,4}(\nu[2]) f_{b,0}(\nu[3])f_{b,1}(\nu[3])\cdots\]
So, the number $\nu[2]$ can be reconstructed as the fixed point of the mapping $f_{b,4}$, and $\nu[3]$ as the fixed point of the mapping $f_{b,0}$. All the other elements of $\nu_w$, including $\nu[1]$, can be computed one-by-one as functions of previously known values. In particular, $\nu[1]=f_{b,3}(\nu[2])$.

For the other fixed point 
\[w=w_b=babab.ab.babab.ab.babab.ab.babab.babab.\ldots=\varphi(w)=\sigma^5(\psi(w)),\]
we also have $w=\sigma^5(\psi(w))$, but since the length of the image of the first symbol $b$ is 5, it means just that
\[w=\psi(\sigma(w)).\]
For the sequence $\nu$, it means that
\[\nu_w=f(\sigma(\nu_w)).\]
So, $\nu[2]$ is the fixed point of $f_{a,1}$, $\nu[3]$ is the fixed point of $f_{b,0}$, and all the other numbers of the sequence $\nu$ can be found starting from them.}
\end{example}

The same idea can be used for every morphism obtained from a ``good'' one by transferring common prefixes of images of symbols to the right: the needed values can be reconstructed from one or several fixed points of mappings $f_{x,i}$.

\subsection{Inseparable types}\label{ss:insep}
%
In this subsection, we propose a method to avoid the situation described in Example \ref{e:aabab}, when the types of elements of the orbit of $w$ are well-defined but cannot be ordered since the relations between words of two given types can be different. We shall show that it happens because of prefixes of these words of bounded length, which can be classified and considered as symbols of a new larger alphabet. The sequence on this new alphabet will inherit all good properties of $w$ and will be separable.

As we have seen above, we may restrict ourselves to a binary order-preserving morphism $\varphi$ such that the last symbols of $\varphi(a)$ and $\varphi(b)$ are different, and the subshift of $\varphi$ is minimal with aperiodic fixed points. Due to Proposition \ref{p:nonprim}, the morphism $\varphi$ is either primitive or of the form $\varphi(a)=axa$, $\varphi(b)=b$ for a finite word $x$ containing $b$.

To discuss the subject, we have to introduce yet another property of morphic subshifts called {\it circularity}. There exist several very close definitions of this property discussed in particular in \cite{ks_circ}; we shall use the following one. The fixed point $w$ of a morphism $\varphi$ (and the whole subshift $L_{\varphi}$) are called {\it circular} if there exists a positive constant $D$ called a {\it synchronization delay} such that in any factor $u$ of $w$ ($L_{\varphi}$) of length at least $D$, there exists a {\it synchronization point}. Here a synchronization point is a place in $u$ where in any occurrence of $u$ to $w$ ($L_{\varphi}$) there is a boundary between images of two symbols: $u=ps$, where $p$ is a suffix of $\varphi(p')$, $s$ is a prefix of $\varphi(s')$, $p's'$ is a factor of $w$ ($L_{\varphi}$).

The smallest value of a synchronization delay can be called {\it the} synchronization delay.

\begin{example}
{\rm
The Thue-Morse word $w_{tm}$ is circular with $D\leq 5$. Indeed, each factor of $w_{tm}$ of length 5 contains one of factors $aa$ or $bb$, and thus a synchronization point between two letters $a$ or two letters $b$. For example, $aabba$ contains two synchronization points, after the first and the third symbol, and appears in $w_{tm}$ only as a suffix of $\varphi_{tm}(bab)$.

The Sierpinski morphism $a \to aba, b \to bbb$ is not circular since for all $n$, the word $b^n$ appears in it and has no synchronization points: the boundaries between images of $b$ in it can pass anywhere. 
}
\end{example}

Note that in our case, when the last letters of images of symbols are all different, a synchronization point $u=ps$ determines a unique decomposition to images of symbols of the whole preceeding prefix $p$ of $u$: we reconstruct it from right to left taking each time the image of symbol ending by the given last letter. At the same time, the suffix $s$, if it is short, may leave some ambiguity if $s$ is the prefix of both images of letters or, in the case of $\varphi(a)$ prefix of $\varphi(b)$ (or vice versa), $s$ is a prefix of $\varphi(ab)$ and of $\varphi(ba)$.

It is well-known that a fixed point of a primitive morphism is circular \cite{mosse,ks_circ}. It is also not difficult to extract from the main result of \cite{ks_rep} and Theorem 12 from \cite{ks_circ} that the non-primitive fixed points from Proposition \ref{p:nonprim} are also circular. So, all morphisms we consider are circular. To be accurate, we redefine the (smallest) synchronization delay $D$ so that, in addition to the main property, each word of length $D$ or more contains both letters: since $w$ is uniformly recurrent, there is no problems with that. We need it to have $|\varphi(u)|\geq D+m-1$ for all $u$ with $|u|\geq D$, where $m$ is the maximum of $|\varphi(a)|, |\varphi(b)|$.

Now let us define another morphism over a greater alphabet preserving all good properties of $\varphi$ and with separable types. To do it, we consider the alphabet $A_D$ of all factors of $w$ of length $D$ and define the trivial isomorphism $\pi: $Fac$_D(w)\mapsto A_D$ which can be naturally extended to $\pi:$Fac$_{D+n}(w) \mapsto A_D^{n+1}$ for all $n$ and to $\pi: L_w \mapsto A_D^{\omega}$, by $\pi(x_1\cdots x_Dy)=\pi(x_1\cdots x_D)\pi(x_2\cdots x_D y)$ for all letters $x_i$ and all words $y$.  Clearly, $\pi$ commutes with the shift $\sigma$ and thus we can consider the subshift $(\pi(L_w),\sigma)$. Moreover, in addition to $\pi^{-1}: A_D \mapsto $Fac$_D(w)$, it is reasonable to consider a simpler mapping $\rho: A_D\mapsto \Sigma_2$, where for each $a \in A_D$, the symbol $\rho(a)$ is its first symbol. The alphabet $A_D$ and words over it inherit the lexicographic order on $\Sigma_2$.

Now, given a morphism $\varphi: \Sigma_2^*\mapsto \Sigma_2^*$, let us define the morphism $\chi: A_D^*\mapsto A_D^*$ as follows: for all $a \in A_D$ such that $a=\pi(u)$ and $\rho(a)=x$ (so that $x$ is the first symbol of $u$), the image $\chi(a)$ is defined as the first $|\varphi(x)|$ symbols of $\pi(\varphi(u))=\pi(\varphi(\pi^{-1}(a)))$. This mapping is well-defined since by the definition of $D$, we have $|\varphi(u)|\geq D+|\varphi(x)|-1$.

\begin{example}
 {\rm
Let us continue Example \ref{e:aabab} (the subshift $L_\varphi$ defined by the morphism $\varphi : a \to aabab$, $b \to 
bba$) and define the respective morphism $\chi$. To do it, we first observe that the synchronization delay of $L_\varphi$ is 5: indeed, the longest word without the synchronization point is $babb$ which can be decomposed both as the factor of $\varphi(ab)$ without the prefix and the suffix of length 2 each, and the factor of $\varphi(bb)$ without the prefix and the suffix of length 1 each. Its continuations $babba$ and $babbb$ disambiguate the situation.

The set Fac$_5(w)$ is of cardinality 17: in the lexicographic order, Fac$_5(L_{\varphi})=\{aaaba,aabab,abaab,ababa,ababb,abbaa,abbab,abbba,$ $baaab,baaba,babaa,$ \\$babba,babbb,bbaaa,bbabb,bbbaa,bbbab\}$. We denote the elements of the alphabet $A_5$, in the same order, as
$A_5=\{a_1,\ldots,a_8,b_1,\ldots,b_{9}\}$: so, $\pi(aaaba)=a_1$ and so on till $\pi(bbbab)=b_{9}$. Note that by the notation, $\rho(a_i)=a$ and $\rho(b_j)=b$ for all well-defined $i$ and $j$.

To define $\chi$, we take the $\varphi$-images of words from Fac$_5(L_{\varphi})$ and then their prefixes of length 9 for words starting from $a$ and of length 7 for words starting from $b$. Then we take $\pi$-images of these words. For example, to find $\chi(a_1)$, we take $\varphi(aaaba)=aababaababaababbbaaabab$, then its prefix of length 9 which is $aababaaba$, then its $\pi$-image $a_2a_4b_3a_3b_2$. So, $\chi(a_1)=a_2a_4b_3a_3b_2$, and continuing the same method, we get
\[\chi:\left\{
\begin{aligned}
        a_1,a_2 &\to a_2a_4b_3a_3b_2,\\
        a_3,a_4,a_5 &\to a_2a_5b_5a_8b_8,\\
        a_6,a_7,a_8 &\to a_2a_5b_5a_8b_{9},\\
        b_1,b_2,b_3,b_4,b_5 &\to b_6b_1a_1,\\
        b_6,b_7 &\to b_7 b_4 a_6,\\
        b_8,b_9 &\to b_7b_4 a_7.
        \end{aligned}\right.\]
Note that for all $x \in A_5$, $\varphi(\rho(x))=\rho(\chi(x))$, so, $\varphi \circ \rho = \rho \circ \chi$. In particular, the two fixed points of $\varphi$, starting from $a$ and from $b$, are $\rho$-images of the two fixed points of $\chi$, starting from $a_2$ and $b_7$. 

The next proposition, following directly from the construction, claims that this is a general situation. } 
\end{example}
\begin{proposition}\label{p:chi}
\begin{enumerate}
 \item 
 If the morphism $\varphi$ is order-preserving on its subshift, then so is $\chi$.
\item
If $\varphi$-images of letters end by different symbols, then different $\chi$-images of symbols of $A_D$ end by different letters of $A_D$. Moreover, these last letters do not occur anywhere else in $\chi$-images of symbols.
\item
There are as many fixed points of $\chi$ as of $\varphi$, and any fixed point of $\varphi$ can be obtained as the $\rho$-image of a fixed point of $\chi$.
\item If $\varphi$-images of letters end by different symbols, then $\rho$ is an isomorphism between $(L_{\chi},\sigma)$ and $(L_{\varphi},\sigma)$, and moreover, it preserves the lexicographic order of words on the respective subshifts.
\item
If $\varphi$-images of letters end by different symbols, then the subshift $(L_{\chi},\sigma)$ is separable.
\end{enumerate}
\end{proposition}
\noindent {\sc Proof.} The first four properties follow directly from the construction. It remains to show the separability of $(L_{\chi},\sigma)$. First let us prove that this subshift is typable. Indeed, for any representation $u=\sigma^{p}(\chi(u'))$, where $u' \in L_{\chi}$ and $0\leq p< |\chi(u'[0])|$, consider the $\rho$-images $v=\rho(u) \in L_{\varphi}$ and $v'=\rho(u') \in L_{\varphi}$. Clearly, $v= \sigma^{p}\varphi(v')$. The morphism $\varphi$ is circular, the last symbols of two images of letters are different, which means that any synchronization point determines all preceeding synchronization points, and so the type of $v$ is well-defined: it is $(v'[0],p)$. But $u=\pi(v)$ and $u'=\pi(v')$ since $v$ and $v'$ belong to $L_{\varphi}$ and by the definition of $\pi$, so, since $v'$ and $p$ are unique, so is $u'$ (and the same $p$). The type of $u$ is thus well-defined as $(u'[0],p)$. Note that we uniquely reconstruct $u'[0]$ even if there are several symbols in $A_D$ with the same $\chi$-image, like in the example above.

It remains to prove that the types in $L_{\chi}$ are comparable. Consider two elements $u_1,u_2 \in L_{\chi}$ of different types $(x_1,p_1)$ and $(x_2,p_2)$. If their first symbols (uniquely defined by types) are different, the order is determined by them. If by contrary the first symbols are the same, let us compare the suffixes of $\chi(x_1)$ and $\chi(x_2)$ which are prefixes of $u_1$ and $u_2$. If $\chi(x_1)\neq \chi(x_2)$, then, by a previous property, the last symbols of these two images are different and do not appear anywhere else in $\chi$-images of letters. So, the order between $u_1$ and $u_2$ is again uniquely determined by the prefixes of $u_1$ and $u_2$ which are suffixes of $\chi(x_1)$ and $\chi(x_2)$. If $\chi(x_1)= \chi(x_2)$ but $p_1\neq p_2$, the same argument holds. At last, if  $\chi(x_1)= \chi(x_2)$, $p_1= p_2$, but $x_1\neq x_2$, we have $u_1<u_2$ if and only if $x_1<x_2$ since the morphism $\chi$ is order-preserving on $L_{\chi}$. This completes the proof of separability of the subshift $(L_{\chi},\sigma)$.
 \hfill $\Box$

\medskip
This construction completes the algorithm allowing to treat every binary pure morphic uniquely ergodic subshift $(L_{\varphi},\sigma)$. First, if the morphism $\varphi$ is not order-preserving, we pass to its square according to Proposition \ref{p:monotoneSq}. Then, if the resulting morphism is not typable because of common suffixes of images of letters, we transfer these common suffixes to to the left as many times as needed due to Proposition \ref{p:suffixes}. At last, we find a synchronization delay and construct the morphism $\chi$ on the extended alphabet which, due to Proposition \ref{p:chi}, has all the desired properties. So, the morphism $f$ constructed as described in the beginning of Section \ref{s:ext} due to Theorem \ref{t:main} gives a equidistributed morphic subshift on the interval $[0,1]$. Again due to propositions from this section, this is the same subshift as it would be for the initial morphism $\varphi$ instead of $\chi$. The value and the infinite sequence corresponding to a fixed 
point of the initial morphism, we can use the method described in Examples \ref{e:aababb} and \ref{e:lettertransf}.

\begin{remark}\label{ex:nonmon}
 {\rm Note that the restriction to the binary alphabet is crucial. On the three-letter alphabet, it is easy to construct a morphism which does not become order-preserving even when we consider its powers: 
$$g: a \to ac, b \to ab, c \to cb.$$
It can be easily seen that $g^n(b)<g^n(a)$ for all $n \geq 1$. Moreover, it is not possible to transfer the common suffix $b$ of $g(b)$ and $g(c)$ to the left since $g(a)$ does not end with $b$. So, we see two technical problems in one example. The extension of the result to general morphisms on larger alphabets is thus an open problem.}
 \end{remark}

\section{Computational tool}\label{s:experiment}
We conclude the paper by a presentation of the software which, given a binary morphism $\varphi$ with an aperiodic uniformly recurrent fixed point, computes the respective morphism $f$ on numbers and first $l$ elements of the sequence corresponding to each of the fixed points of $\varphi$. The code is available at 

\noindent  {\tt
https://www.i2m.univ-amu.fr/perso/anna.frid/MorphismsOnReals/mp.py}.

\normalsize
\noindent A web page where the computation can be done online for a relatively small input (images of letters not longer than about 5 letters, if the morphism is not order-preserving, or about 25 letters, if it is) is 

\noindent {\tt https://www.i2m.univ-amu.fr/perso/anna.frid/MorphismsOnReals/mp.html}.

1) Given a binary morphism $\varphi$, we first check if it is primitive: in the binary case, it is clearly sufficient to check if the square of its matrix is positive. If the morphism is not primitive, we continue to consider it if it falls into the case of Proposition \ref{p:nonprim}. Then we check if $\varphi$ admits a fixed point starting with each letter, that is, if $\varphi(x)$ starts with $x$ for some letter $x$.  
To check fixed points for aperiodicity, we use the result of \cite{seebold}
and eliminate periodic fixed points corresponding $\varphi(a) = a(ba)^m$ and $\varphi(b) = b(ab)^n$ for some $m,n \geq 0$ such that $m+n \geq 1$, and those whose  images are powers of a common word of length at least $2$. Among the  non-primitive morphisms from Proposition \ref{p:nonprim}, we eliminate those of the form $\varphi(a)=a(b^ma)^n$ and $\varphi(b)=b$ (and of course the symmetric case). 

2) In the primitive case, we check if the morphism is order-preserving. If one image of a letter is not a prefix of the other one, the check is straighforward; if it is, we can transfer the common prefix to the end of both images until they can be directly compared. Corollary \ref{c:opres} assures that the property of being order-preserving is stable under this operation. If the morphism is not order-preserving, we pass to its square due to Proposition \ref{p:monotoneSq}, even though it considerably increases the complexity of the computation. Then we check if the two images of letters have a common suffix, at if it is the case, we transfer it to the beginning of the images as it was described in Proposition \ref{p:pssp}. Due to further results of Subsection \ref{ss:cs}, it is sufficient to repeat this procedure a finite number of times, and the resulting morphism $\psi$ remains order-preserving.

3) The morphism $\mu$ considered at this stage (here $\mu$ can be the initial $\varphi$, or its square, or the $\psi$ obtained from $\varphi$ or $\varphi^2$ by transferring common suffixes to the left) is circular. 
Now we have to compute a synchronization delay $D$ of $\mu$ to check if $\mu$ is separable; if it is not, we will have to use yet another morphism on a greater alphabet. This is a slow part of the computation, especially as in the general case, there is no known upper bound for $D$. We only know from a recent paper by Klouda and Medkov\'a \cite{km} that for a uniform binary morphism, $D$ is bounded by $m^3$, where $m$ is the morphism length. So, in the general case, we unfortunately do not have an upper bound for the complexity of the following procedure.

First, we find all factors of length 2 of the subshift as follows: starting from the set $A$ of all factors of length 2 of images of letters, we expand $A$ while there are new words of length 2 situated at the boundary between images of letters in words $\mu(a)$, $a\in A$. Clearly, as soon as there are no new words of length 2 obtained like that, the set of factors of length 2 of the subshift $L$ is complete.

Now, starting with the set of factors of length 2, we use the following fact: If $m$ is the shorter length of an image of a letter, then every factor of $L$ of length $ml+1$ is contained in a $\mu$-image of a factor of $L$ of length $l+1$. If $m\geq 2$, this fact is sufficient to find all factors of $L$ together with their types for any given length greater than 2. If $m=1$ and $\mu$ is primitive, we can pass to its square to get $m>1$. At last, for the non-primitive case of $\varphi(a)=axa$, $\varphi(b)=b$, we need the following two lemmas which allow to initialize and continue the process of finding all factors of any given length. 

\begin{lemma} Let $\varphi$ be a binary morphism such that  $\varphi(a)=axa$ and $\varphi(b)=b$, where $x$ is a finite word containing $b$. If we denote by $h_b$ the highest power of $b$ appearing in $x$, then the factors of length $h_b+2$ of $L_{\varphi}$ are exactly those of  $\varphi^2(a)$.
\end{lemma}

\noindent {\sc Proof.} 
Any factor of $\varphi^2(a)$ is also a factor of $L_{\varphi}$. Conversely, let $u$ be a factor of length $h_b+2$ of the subshift. Then $u$ contains an occurrence of $a$ and thus any its given occurrence has common letters with an occurrence of $\varphi(a)$. If it is contained in $\varphi(a)$, the statement is proved; if not, since $|\varphi(a)| \geq h_b+2$, the word $u$ overruns $\varphi(a)$ from at most one side, so, it is a factor of $\varphi(a)b^k \varphi(a)$ for some $k \leq h_p$. In particular, $\varphi(a)b^k \varphi(a)$ is a factor of the subshift, and for that, its preimage $ab^k a$ had to appear, at some point, in $\varphi(a)$. 
So, $\varphi(a)b^k \varphi(a)$ and its factor $u$ are factors of $\varphi (\varphi(a))$. \hfill $\Box$ 

\begin{lemma}
 Let $\varphi$ be a binary morphism such that  $\varphi(a)=axa$ and $\varphi(b)=b$, where $x$ is a finite word containing $b$. If we denote by $h_b$ the highest power of $b$ appearing in $x$, then for every $l\geq h_b+2$, every factor of $L_{\varphi}$ of length less than $p(l)$ is a factor of some $\varphi(y)$, where $y$ is a factor of $L_{\varphi}$ of length at most $l$. Furthermore, we have $p(l)>l$. Here $p(l)$ is defined by
 $$p(l)  =  1 + q(h_b+|\varphi(a)|) + r,$$
  where $q$ and $r$ stand for the quotient and the rest in the euclidean division 
of $l-1$ by  $h_b+1$.
\end{lemma}
\noindent {\sc Proof.} 
The shortest word that can be written by concatenating $l-1$ $\varphi$-images of letters and which does not contain the word $b^{h_b+1}$ is $(\varphi(b)^{h_b}\varphi(a))^q\varphi(b)^r$ of length $q(h_b+\vert \varphi(a) \vert)+r$. As a consequence, the length of the shortest factors of $L_{\varphi}$ that lie on $l+1$ (and not less) $\varphi$-images of letters is at least $q(h_b+\vert \varphi(a) \vert)+r+2=p(l)+1$, meaning that all factors of length $p(l)$ or less lie on at most $l$ $\varphi$-images of letters.  It is easy to check that for $l\geq h_b+2$ we have $p(l)>l$. \hfill $\Box$

\medskip
So, in each situation, we can find all factors of any given length together with their types. A synchronization delay is reached as soon as every word appears in the list with only one type. It is reasonable to check all the shorter lengths and find {\it the} smallest synchronization delay $D$ not an upper bound for it, to consider a smaller alphabet and to have a nicer output (and probably to gain in computation time). 

4) We have obtained the set of factors of length $D$, and each of them corresponds to a type. However, the same type may correspond to several words, and if they are not lexicographically consecutive, the morphism is not separable. If it is the case, we pass to a morphism $\chi$ on a larger alphabet as described in Subsection \ref{ss:insep}.

5) The morphism considered at this stage satisfies the conditions of Theorem \ref{t:main} and so we construct the morphism on intervals as described in Section \ref{s:ext}.

6) To find numeric sequences corresponding to fixed points of the initial morphism $\varphi$, we start with finding fixed points of respective mappings $f_{c,p}$ as described in Examples \ref{e:tm}, 
\ref{e:aababb}, \ref{e:lettertransf}.

At last, note that the Perron-Frobenius eigenvalue is an algebraic number, making it possible for a mathematical software to do exact computations at each stage, and to print the outputs with arbitrary precision.

As we have discussed, the algorithm we use is not very fast. First, there is no general upper bound for the synchronization delay $D$, and at the same time, computing $D$ is the slowest part of the process. To avoid it, it would be nice to invent a faster way to check the subshift for separability. It would be also helpful to learn how to deal directly with order-reversing morphisms, since taking the morphism square before looking for separability slows down this slowest part of computation. We leave these questions to further research.

\end{document}